\documentclass{article}

     \usepackage[nonatbib]{neurips_2020}

\usepackage[utf8]{inputenc} 
\usepackage[T1]{fontenc}    
\usepackage{hyperref}       
\usepackage{url}            
\usepackage{booktabs}       
\usepackage{amsfonts}       
\usepackage{nicefrac}       
\usepackage{microtype}      
\usepackage{graphics,graphicx,epsfig,algorithm,algpseudocode}
\usepackage{latexsym,amsmath,amsfonts,amssymb,color}
\usepackage{graphicx}
\newcommand\tabcaption{\def\@captype{table}\caption}
\newcommand{\be}{\begin{eqnarray*}}
\newcommand{\ee}{\end{eqnarray*}}
\newcommand{\ibe}{\begin{eqnarray}}
\newcommand{\iee}{\end{eqnarray}}

\title{Fast Numerical Simulation of Allen--Cahn Equation}

\author{%
  Yongho Kim \\
  Institute of Computer Science\\ University of Hildesheim\\ Samelsonplatz 1, 31141 Hildesheim, Germany\\
  \texttt{kimy@uni-hildesheim.de} \\
}

\begin{document}

\maketitle

\begin{abstract}
Simulation speed depends on code structures, hence it is crucial how to build a fast algorithm. We solve the Allen--Cahn equation by an explicit finite difference method, so it requires grid calculations implemented by many for-loops in the simulation code. In terms of programming, many for-loops make the simulation speed slow. To solve the problem, we propose a model architecture containing a pad and a convolution operation for the Allen--Cahn equation. Also, the GPU operation is used to boost up the speed more. In this way, the simulation of other differential equations can be improved. In this paper, various numerical simulations are conducted to confirm that the Allen--Cahn equation follows motion by mean curvature and phase separation in two-dimensional and three-dimensional spaces. Finally, we demonstrate that our algorithm is much faster than an unoptimized code and the CPU operation.
   
\end{abstract}

\section{Introduction}

The Allen--Cahn (AC) equation is a reaction-diffusion equation composed of the reaction term $-{F'(\phi(\textbf{x},t))}/{\epsilon^2}$ and the diffusion term $\Delta \phi(\textbf{x},t)$:
\begin{equation}
\phi_t(\textbf{x},t)=-\frac{F'(\phi(\textbf{x},t))}{\epsilon^2}+\Delta \phi(\textbf{x},t), \quad \textbf{x} \in \Omega, t>0, \label{oriACeq}
\end{equation}
where $\phi(\textbf{x},t)$ is order parameter which is defined as the difference in concentration of the two components in a mixture and $F(\phi)$ is double well potential energy function with minimum values at $-1$ and $1$, and its form is $F(\phi)=0.25(\phi^2 -1)^2$. $\epsilon$ is thickness of transition layer and which is small positive constant value.  The AC equation is first introduced in a research on the phase separation of binary iron alloys \cite{AC1979}. The AC equation has studied and applied to various fields such as image inpainting \cite{YDJSJ2015,AA2018,ZJJ2008}, image segmentation \cite{MVK2004,DA2009}, crystal growth \cite{JCX2019,XY2019,XQ2020}. \\

With the development of computer hardware such as GPU(Graphics Processing Unit) and memory cards, neural networks are applied in a wide range of research areas such as computer vision, natural language processing, and numerical analysis. As GPU operations outperform CPU(Central Processing Unit) performance in multi-tasks and high-dimensional problems, open source machine learning libraries such as Pytorch provide a variety of neural networks using GPU and are useful to build the architecture combined with neural networks and numerical methods. Thus, many researches use machine learning libraries. M. Raissi et al.\cite{r1} proposed physics informed neural networks combined by multi-layer perceptrons and numerical methods to solve nonlinear partial differential equations, S. Karumuri et al.\cite{r2}  introduced a solver-free approach for Stochastic partial differential equations, and L. Yang et al.\cite{r3} proposed a Bayesian physics informed neural network. In this paper, we propose a structure using padding and convolution operation for the GPU calculation of the AC equation and demonstrate the validity of the proposed structure by verifying the result with Python code that has the same mathematically meaning.\\

This paper is organized as follows. In Section \ref{NuSol}, we present an explicit finite difference method to solve the AC equation which implemented by CPU and GPU algorithms. In Section \ref{Nuexp}, the numerical simulations including a motion by mean curvature, phase separation, and temporal evolutions of various initial shapes are introduced as well as the runtime results between CPU and GPU operations are compared. Finally, conclusions are drawn in Section \ref{con}.

\section{Numerical solutions} \label{NuSol}
In this section, we present an explicit finite difference method to solve the AC Eq. \eqref{oriACeq}. Also, we give an explanation of each algorithm for CPU and GPU computing. For simplicity of expression, we describe a numerical scheme for the AC equation in two dimensions (2D), and the definition in three dimensional (3D) space can be easily extended and considered. A computational domain is defined using a uniform grid of size $h=1/N_x$ and $\Omega_h=\{ (x_i, y_j)=(a+(i-0.5)h, c+(j-0.5)h) \}$ for $1 \leq i \leq N_x$, $1 \leq j \leq N_y$ is the set of cell centers. Here, $N_x, N_y$ are mesh size on computational domain $(a,b)\times(c,d)$. For the definition of the boundary condition, we define the extended computational domain as follows:
\begin{equation*}
\Omega_h=\{ (x_i, y_j)=(a+(i-1.5)h, c+(j-1.5)h) \},
\end{equation*}
for $1 \leq i \leq N_x+2$ and $1 \leq j \leq N_y+2$. Let $\phi_{ij}^n$ be approximations of $\phi(x_i,y_j,n\Delta t)$, where $\Delta t=0.1h^2$ is temporal step size, $T$ is a final time, and $N_t$ is a total number of time steps. The boundary condition is zero Neumann boundary condition: 
\begin{equation*}
\begin{aligned}
&\phi_{i,1}^n=\phi_{i,2}^n, \quad \phi_{i,N_y+2}^n=\phi_{i,N_y+1}^n \qquad 1 \leq i \leq N_x,\\
&\phi_{1,j}^n=\phi_{2,j}^n, \quad \phi_{N_x+2,j}^n=\phi_{N_x+1,j}^n \qquad 1 \leq j \leq N_y.
\end{aligned}
\end{equation*}
We define the thickness of transition layer $\epsilon$ in Eq. \eqref{oriACeq} as $\epsilon_m$ \cite{JCHLDJJK2009}:
\begin{equation}
\epsilon_m = \frac{hm}{2\sqrt{2} \tanh^{-1}(0.9)}, \label{emvalue}
\end{equation}
where $m$ is the number of grids representing the thickness.

First of all, we show the error results obtained by Eq. \eqref{ERReq} to check the difference between the numerical results of CPU:Python (baseline) and GPU:Pytorch. In Table \ref{tab:3},  all the errors for any cases are less than $1.0e$-$6$.

\begin{table}[htbp]
 \begin{center}
   \caption{Errors of various numerical simulations with the baseline and ours}\label{tab:3}
   \resizebox{\textwidth}{!}{
   \begin{tabular}{c c c c c c c}
   \hline
       &  \multicolumn{6}{c}{Initial value} \\
     \cline{2-7}
     $Dimension$ & $separation$ &  $dumbbell$ & $ circle/sphere$ & $maze$ & $star$ & $torus$\\ 
     \hline
     2D&  $1.75\times 10^{-6}$ & $7.03\times 10^{-7}$& $5.51\times 10^{-7}$& $2.22\times 10^{-7}$& $6.55\times 10^{-7}$& $7.52\times 10^{-7}$\\
     3D&  $3.01\times 10^{-6}$ & $1.20\times 10^{-6}$& $1.11\times 10^{-6}$& $3.36\times 10^{-6}$& $1.56\times 10^{-6}$& $1.91\times 10^{-6}$\\
     \hline
   \end{tabular}}
   \vspace{1ex}
 \end{center}
\end{table} 

To estimate the error of the CPU and GPU codes, we use the following defined $Err$
\begin{equation}
Err = \frac{1}{n}\sum_{t=1}^{n}{\sqrt{avg((\phi_{gpu}^{t}-\phi_{cpu}^{t})^2)}} \label{ERReq}
\end{equation}
where $avg(X)$ is the average of elements in an array $X$ and $t$ is the time step. As the result, it is confirmed that there is little difference between the two algorithms.

\subsection{Numerical solutions on CPU (baseline)} 

The AC equation \eqref{oriACeq} is discretized using the explicit finite difference method as:
\begin{equation}
\begin{aligned}
\frac{\phi_{ij}^{n+1}-\phi_{ij}^{n}}{\Delta t}&=\frac{\phi_{ij}^{n}-(\phi_{ij}^{n})^3}{\epsilon^2}+\Delta_h\phi_{ij}^n\\
\phi_{ij}^{n+1}&=\phi_{ij}^{n}+\Delta t
\left(\frac{\phi_{ij}^{n}-(\phi_{ij}^{n})^3}{\epsilon^2}+\Delta_h\phi_{ij}^n\right)\\
\phi_{ij}^{n+1}&=(1+\alpha)\phi_{ij}^n - \alpha(\phi_{ij}^n)^3+\Delta t\Delta_h\phi_{ij}^n,
\end{aligned} \label{pytorchACeq}
\end{equation}
where $\alpha = {\Delta t}/{\epsilon^2}$ and $\Delta_h \phi_{ij}^n=({ \phi^n_{i-1,j}+\phi^n_{i+1,j}+\phi^n_{i,j-1}+\phi^n_{i,j+1}-4\phi^n_{i,j}})/{h^2}$. The baseline algorithm is implemented in Eq. \eqref{pytorchACeq} using Numpy (CPU array).

\subsection{Numerical solutions on GPU (Pytorch)} 
\begin{figure}[t!]
\centering
\includegraphics[width=0.7\columnwidth]{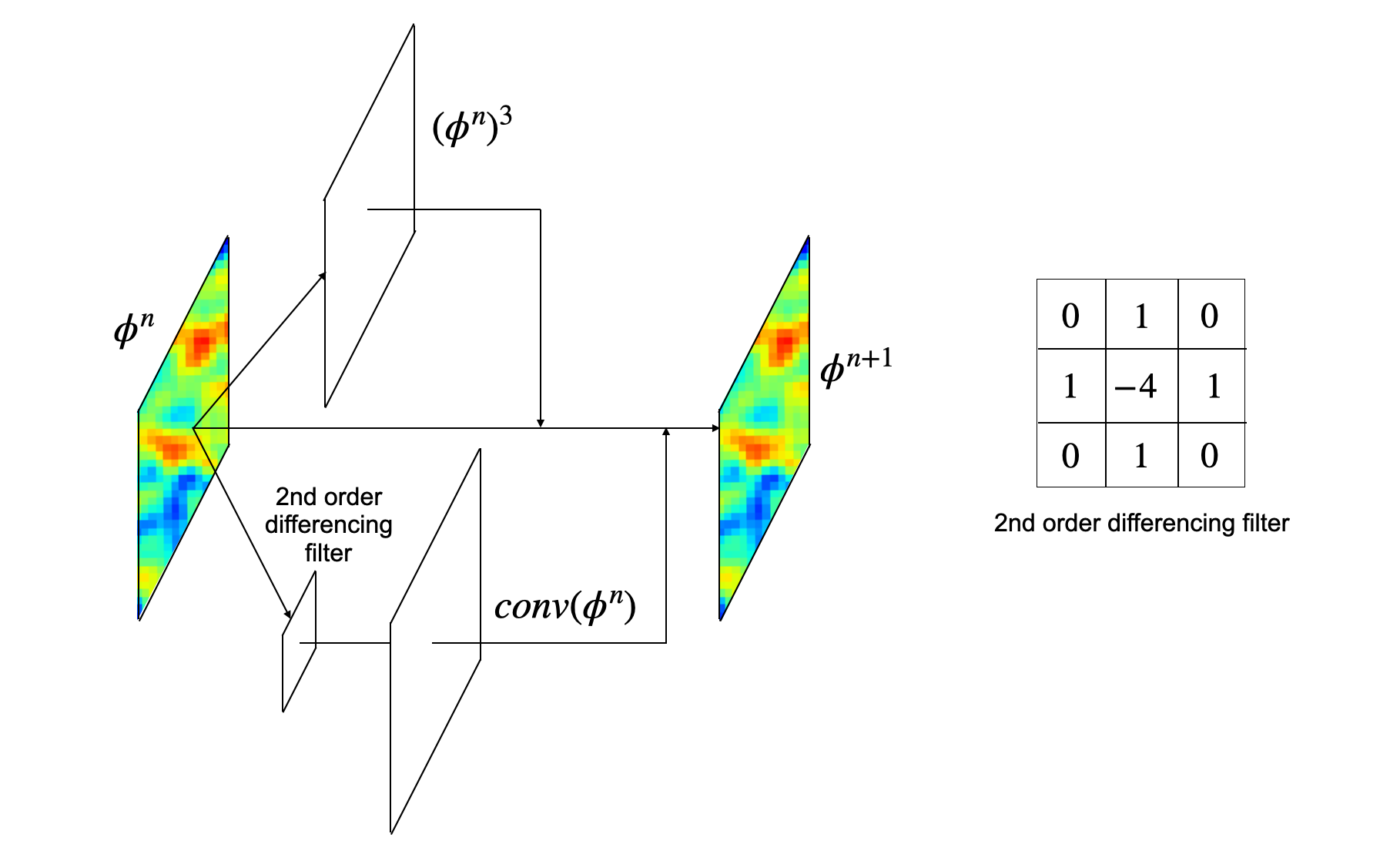}
\caption{Schematic of the process to obtain $f$}
\label{fig_pmat}
\end{figure}

For GPU computing, the AC equation \eqref{oriACeq} can be expressed using Pytorch and the algorithm can be represented as Fig. \ref{fig_pmat}.

\begin{equation}
\begin{aligned}
\frac{\phi^{n+1}-\phi^{n}}{\Delta t}&=\frac{\phi^{n}-(\phi^{n})^3}{\epsilon^2}+\Delta_h\phi^n\\
\phi^{n+1}&=\phi^{n}+\Delta t
\left(\frac{\phi^{n}-(\phi^{n})^3}{\epsilon^2}+\Delta_h\phi^n\right)\\
\phi^{n+1}&=(1+\alpha)\phi^{n}-\alpha(\phi^{n})^3+conv(\phi^n)\quad
\left(\alpha=\frac{\Delta t}{\epsilon^2},\, conv(\phi^n)=\Delta t\Delta_h\phi^n\right)\\
\phi^{n+1}&=f(\phi^n) \label{GPUcodeNet}
\end{aligned}
\end{equation}
where $f$ is a Pytorch model. The main algorithm is a set of steps in Eq. \eqref{GPUcodeNet} using Pytorch. In this algorithm, \textit{nn.ReplicationPad2d} (or \textit{nn.ReplicationPad3d}) is applied to satisfy the boundary condition by padding the $\phi^n$ input using replication of the input boundary. Also, the convolutional operator \textit{F.conv2d} (or \textit{F.conv3d}) with the 2nd order differencing filter is used to calculate the diffusion term $\Delta\phi^n$.
The model can choose  an operation mode between GPU and CPU by \textit{to(device)}. If device=cuda:0, the model is implemented on GPU, otherwise on CPU. All the codes are available from the first author's GitHub webpage: https://github.com/kimy-de/gpuallencahn

\section{Numerical experiments} \label{Nuexp}
In this section, we perform the following numerical tests in 2D and 3D: $\epsilon_m$ effect, the motion by mean curvature effect with various initial shapes [circle (sphere in 3D), dumbbell, star, torus, maze] and phase separation.  In the $\epsilon_m$ effect test, we compare the numerical solutions to analytic values to find a proper $\epsilon_m$ value. We simulate a phenomenon that follows motion by mean curvature in various initial shapes and phase separation with random initial condition. We perform the simulation on the following specifications: Intel(R) Core(TM) i9-9900K CPU @3.60GHz, 32GB RAM / NVIDIA GeForce RTX 2080 Super.

\subsection{Initial conditions} 
We consider the 2D and 3D initial conditions introduced in this section: To check the $\epsilon_m$ effect, we measure the circle's (sphere in 3D) radius that changes with time and take the initial conditions:
\begin{subequations}
\begin{align}
\phi(x,y,0)&=\tanh \left( \frac{R_0-\sqrt{(x-0.5)^2+(y-0.5)^2}}{\sqrt{2}\epsilon} \right), \qquad \textrm{in 2D} \label{2Dcircleini}\\ 
\phi(x,y,z,0)&=\tanh \left( \frac{R_0-\sqrt{(x-0.5)^2+(y-0.5)^2+(z-0.5)^2}}{\sqrt{2}\epsilon} \right), \qquad \textrm{in 3D} \label{3Dsphereini}
\end{align}
\end{subequations}
where $R_0$ is the initial radius of a circle (sphere in 3D).

In various tests to check the motion by the mean curvature, the initial conditions for the circle and sphere refer to Eqs. \eqref{2Dcircleini} and \eqref{3Dsphereini}, respectively. A dumbbell is constructed by following the initial conditions (o/w: otherwise)
\begin{subequations}
\begin{align}
\phi(x,y,0)&=
\begin{cases}
1.0, \quad \textrm{if}~~(0.4<x<1.6) ~~\textrm{and}~~ (0.4<y<0.6) \\
1+\tanh \left( \frac{R_0-\sqrt{(x-0.3)^2+Y}}{\sqrt{2}\epsilon} \right)+\tanh \left( \frac{R_0-\sqrt{(x-1.7)^2+Y}}{\sqrt{2}\epsilon} \right),~\textrm{o/w},
\end{cases}  \textrm{in 2D} \label{2Ddumbellini}\\ 
\phi(x,y,z,0)&=
\begin{cases}
1.0, \quad \textrm{if}~~(0.4<x<1.6) ~~\textrm{and}~~ (0.4<y,z<0.6) \\
1+\tanh \left( \frac{R_0-\sqrt{(x-0.3)^2+YZ}}{\sqrt{2}\epsilon} \right)+\tanh \left( \frac{R_0-\sqrt{(x-1.7)^2+YZ}}{\sqrt{2}\epsilon} \right),\textrm{o/w},
\end{cases} \textrm{in 3D}  \label{3Ddumbellini}
\end{align}
\end{subequations}
where $R_0$ is initial radius of both side of dumbbell's circle (sphere in 3D) and for simplicity of expression, $Y=(y-0.5)^2$ and $YZ = (y-0.5)^2+(z-0.5)^2$.

The initial conditions of a star shape are defined as:
\begin{subequations}
\begin{align}
\phi(x,y,0)&=\tanh \left( \frac{0.25+0.1\cos(6\theta)-\sqrt{(x-0.5)^2+(y-0.5)^2}}{\sqrt{2}\epsilon} \right), \qquad \textrm{in 2D} \label{2Dstarini}\\ 
\phi(x,y,z,0)&=\tanh \left( \frac{0.7+0.2\cos(6\theta)-\sqrt{x^2+y^2+z^2}}{\sqrt{2}\epsilon} \right), \qquad \textrm{in 3D} \label{3Dstarini}
\end{align}
\end{subequations}
where 
\begin{equation*}
\theta=
\begin{cases}
\tan^{-1} \left( \frac{y-0.5}{x-0.5} \right), \quad \textrm{if}~~(x>0.5) \\
\pi + \tan^{-1} \left( \frac{y-0.5}{x-0.5} \right), \quad \textrm{o/w}.
\end{cases} ~\textrm{in 2D}, \qquad
\theta=
\begin{cases}
\tan^{-1} \left( \frac{z}{x} \right), \quad \textrm{if}~~(x>0.5) \\
\pi + \tan^{-1} \left( \frac{z}{x} \right), \, \textrm{o/w}.
\end{cases} \textrm{in 3D}
\end{equation*}
and we use different domain sizes in 2D and 3D, so the center of the star depends on the dimensions.

And a torus shape is given by
\begin{subequations}
\begin{align}
\phi(x,y,0)&=-1+\tanh \left( \frac{R_1-\sqrt{XY}}{\sqrt{2}\epsilon} \right) - \tanh \left( \frac{R_2-\sqrt{XY}}{\sqrt{2}\epsilon} \right), \qquad \textrm{in 2D} \label{2Dtorusini}\\
\phi(x,y,z,0)&=\sqrt{z^2+\left(\sqrt{x^2+y^2}-R_1\right)^2}-R_2, \qquad \textrm{in 3D} \label{3Dtorusini}
\end{align}
\end{subequations}
where $R_1$ and $R_2$ are the radius of major (outside) and minor (inside) circles, respectively. And, for simplicity of expression, $XY=(x-0.5)^2+(y-0.5)^2$.

The last initial conditions of a maze shape is complicated to describe its equation, so refer to the code in Appendix.

Random initial conditions for confirming the phase separation of the AC equation are
\begin{subequations}
\begin{align}
\phi(x,y,0)&=0.1\textrm{rand}(x,y), \qquad \textrm{in 2D} \label{2Dtrnadini}\\
\phi(x,y,z,0)&=0.1\textrm{rand}(x,y,z), \quad ~\textrm{in 3D} \label{3Drandini}
\end{align}
\end{subequations}
here the function rand$(x,y)$ has a random value between $-1$ and $1$.

\subsection{Simulations in 2 dimensional space} 
Unless otherwise stated, we use the following parameters: mesh size $N_x=N_y=200$, space step size $h=1/N_x$, time step size $\Delta t = 0.1h^2$, and computational domain $\Omega = (0,1)\times(0,1)$. We should find a proper thickness of transition layer as defined Eq. \eqref{emvalue}. First, we find an appropriate $\epsilon_m$ by comparing the numerical solution with the exact solution for the radius that decreases due to the motion by mean curvature in the circle initial condition in Eq. \eqref{2Dcircleini}.

\begin{figure}[t!]
\centering
\includegraphics[width=0.8\columnwidth]{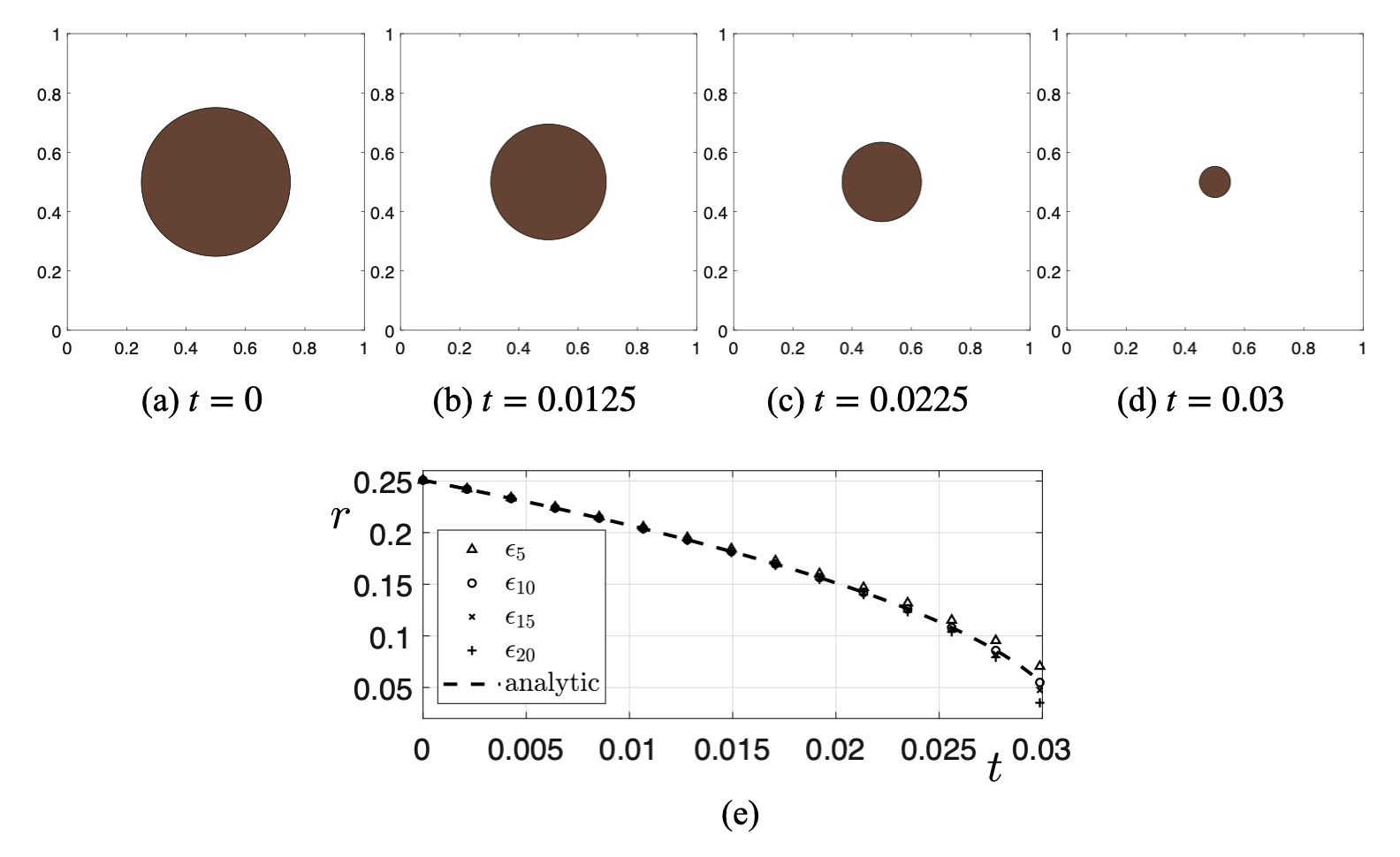}
\caption{(a)-(d) The changes in the circle over time are shown, and each time is described below in the figures. (e) Change of radius for various $\epsilon_m$.} \label{fig_circle_results}
\end{figure}

\begin{figure}[t!]
\centering
\includegraphics[width=1\columnwidth]{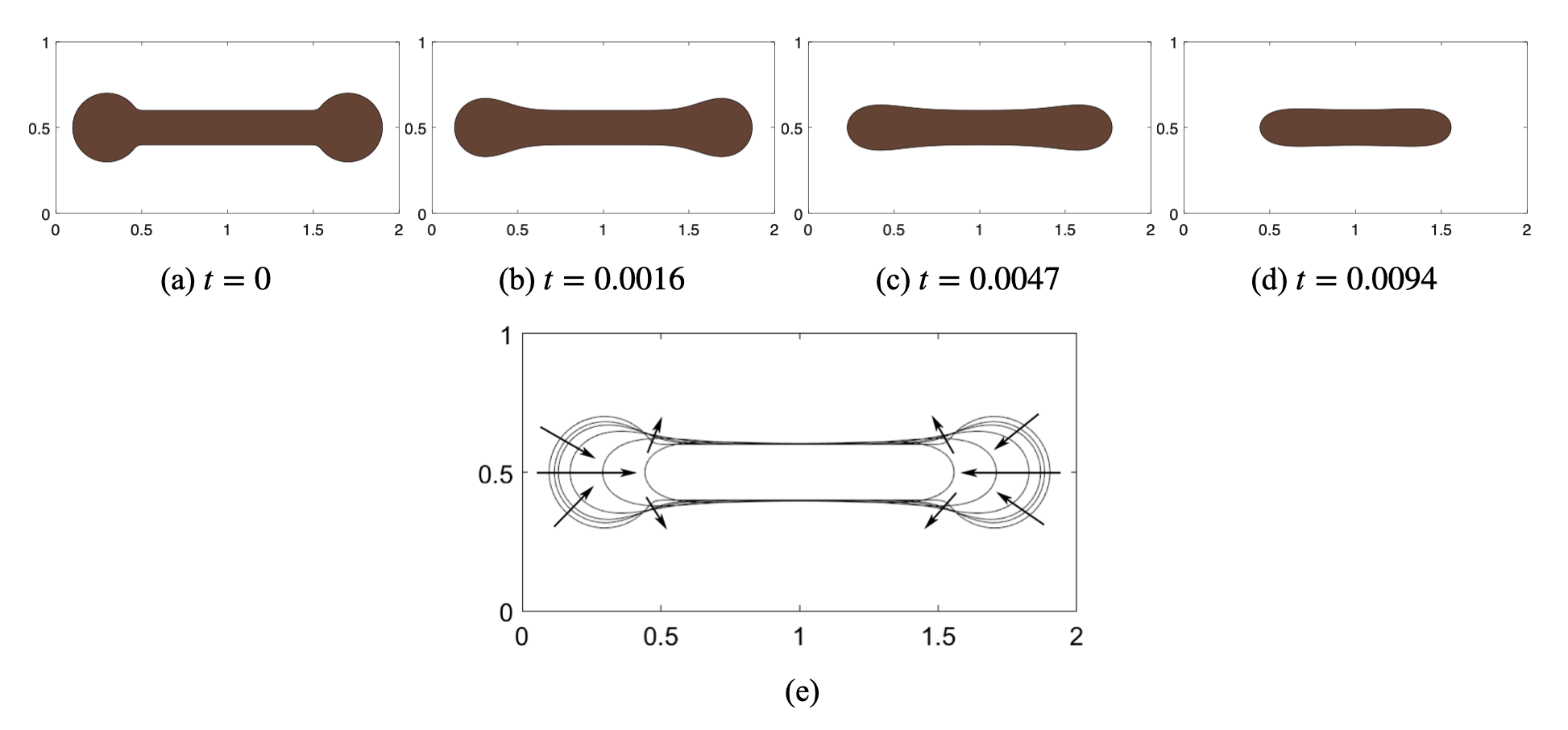}
\caption{(a)-(d) Time evolution of dumbbell shape. And each time is described below in the figures.  (e) Contour lines over time are shown by overlapping. It can be seen that the changing by the mean curvature flow.} \label{fig_dumbell_results} 
\end{figure}

\begin{figure}[htbp!]
\centering
\includegraphics[width=0.82\columnwidth]{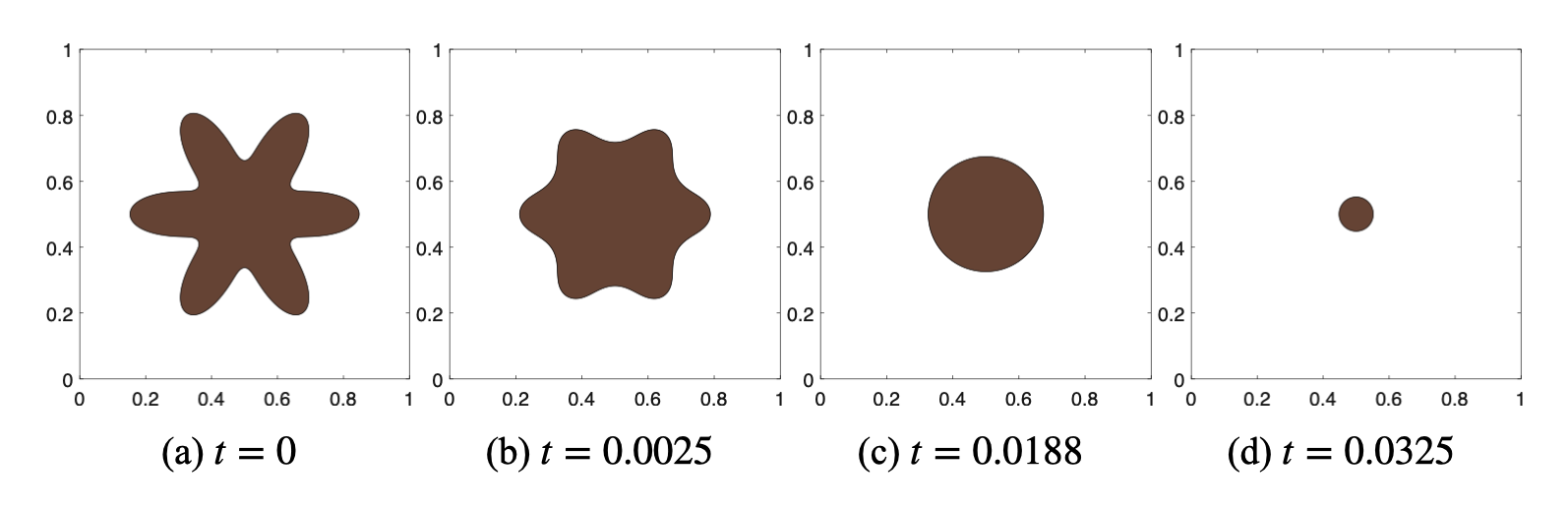}
\caption{Time evolution of a star shape. As in the shape of a dumbbell, it can be seen that it changes with the mean curvature flow.} \label{fig_star_results} 
\end{figure}

\begin{figure}[htbp!]
\centering
\includegraphics[width=0.82\columnwidth]{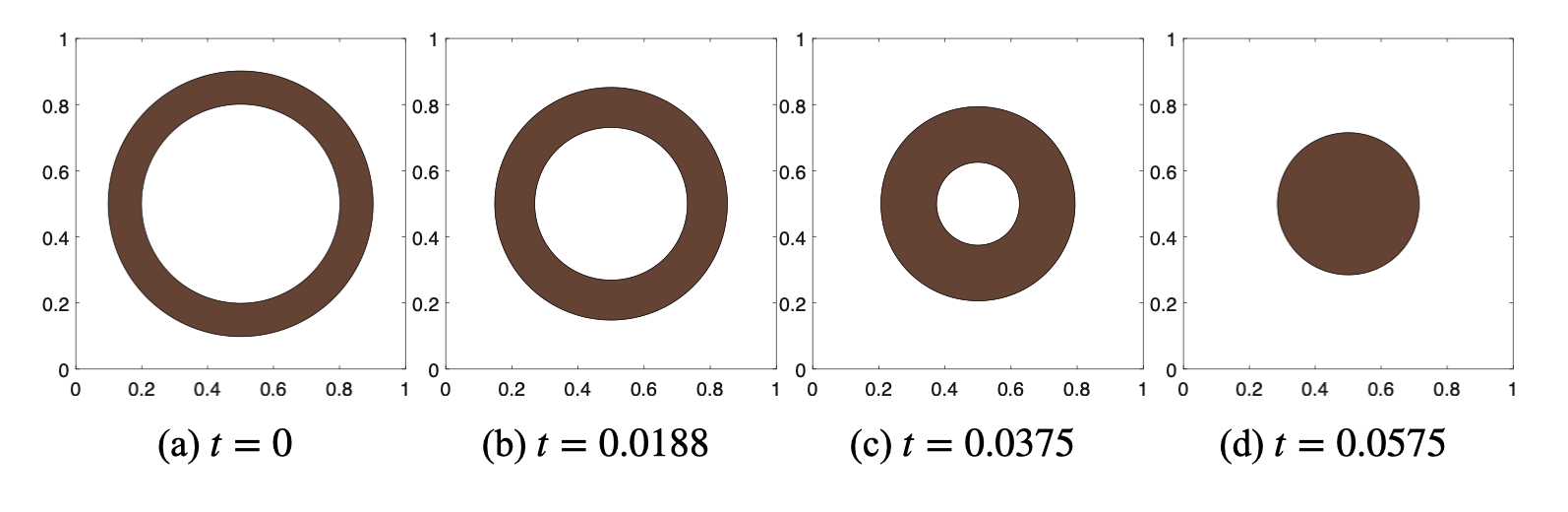}
\caption{Time evolution of a torus shape. And each time is described below in the figures. Inner circle has larger curvature than outer circle, so inner circles shrinks faster than outer circle.} \label{fig_torus_results} 
\end{figure}

\begin{figure}[t!]
\centering
\includegraphics[width=0.82\columnwidth]{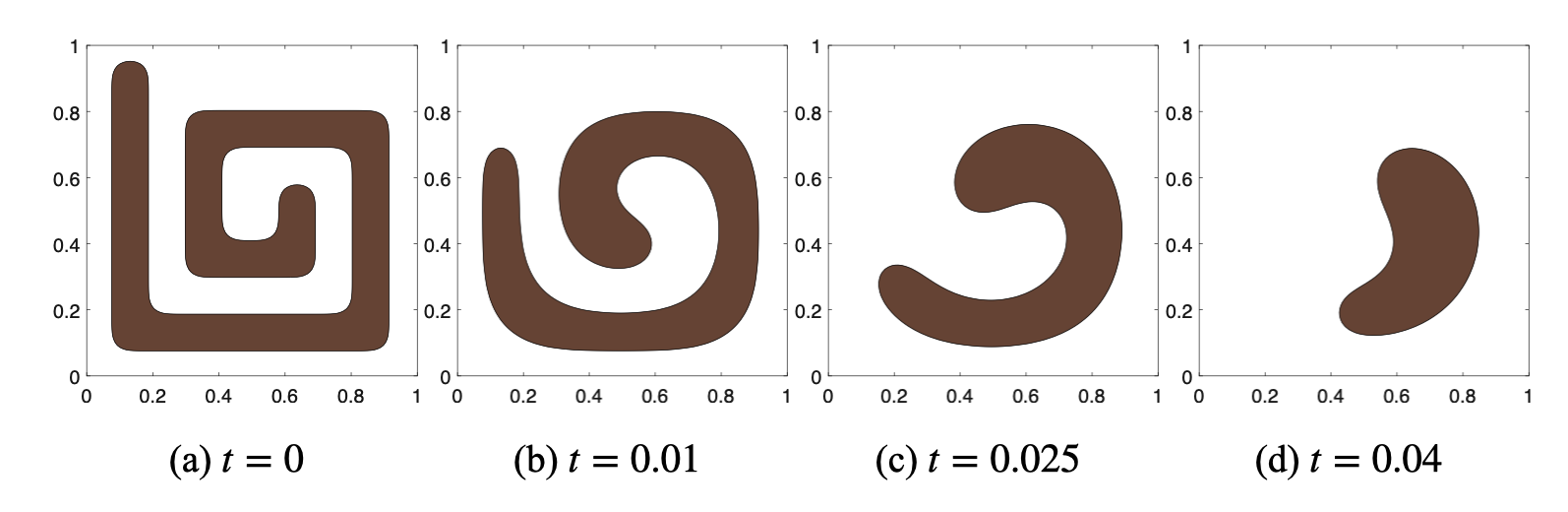}
\caption{Time evolution of a maze shape on $100 \times 100$ mesh size with $\epsilon_5$, and $T=0.04$.} \label{fig_complex_results} 
\end{figure}

\begin{figure}[t!]
\centering
\includegraphics[width=0.82\columnwidth]{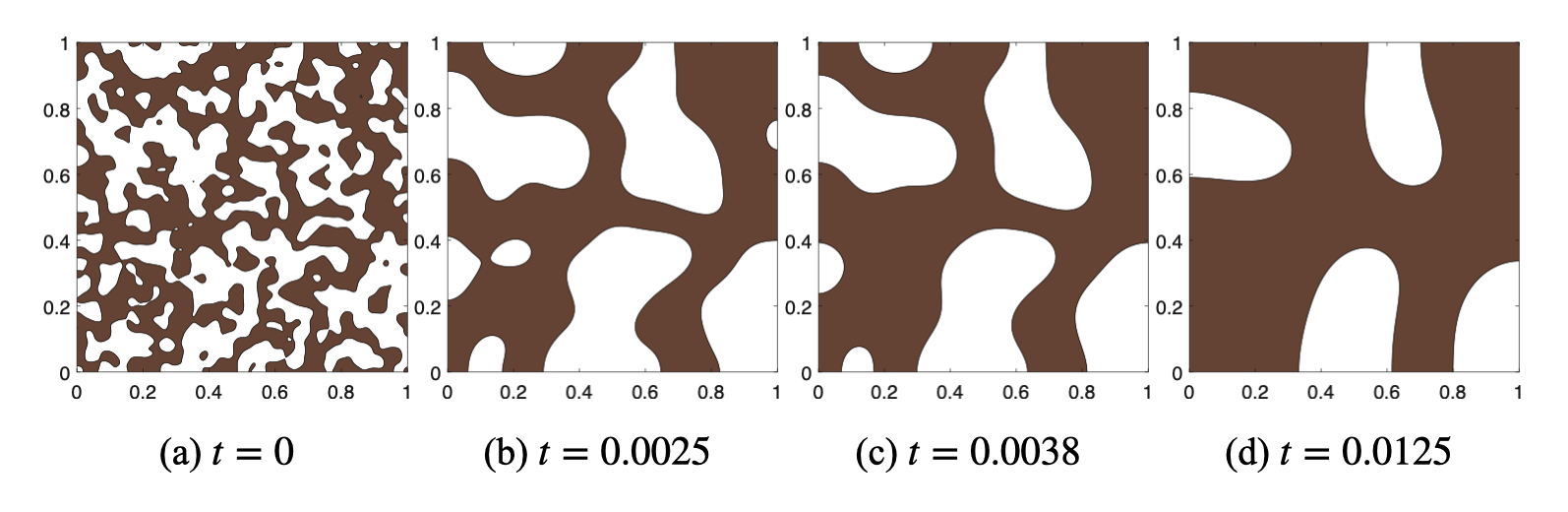}
\caption{Time evolution of phase separation with a random initial condition.} \label{fig_separation_results} 
\end{figure}

Figure \ref{fig_circle_results}(a)-(d) show a circle shrinking with motion by mean curvature based on $\epsilon_{10}$. The exact solution of the radius $r$ decreasing with time evolution can be calculated \cite{YHDJ2010}
\begin{equation}
r(t)=\sqrt{r_0^2+2(1-d)t}, \label{exactradiusEQ}
\end{equation}
where $r_0$ is the initial radius of the circle and $d$ is a dimension, and $t$ is time. We simulate various $\epsilon_m$ until the final time $T=0.03$. As shown in Fig. \ref{fig_circle_results}, when $\epsilon_{10}$ is used, it is found that the exact solution and the numerical solution are most similar in the experiment. Therefore, the other tests are performed using $\epsilon_{10}$ except for the case of changing the grid in 2D.
Figure \ref{fig_dumbell_results} shows the temporal evolution of the dumbbell shape and the initial condition is shown in Eq. \eqref{2Ddumbellini}. For resolution, we use $N_x=400$ and $N_y=200$ on the computational domain $\Omega=(0,2)\times(0,1)$. The final time $T$ of the simulation is $0.0094$ and $R_0=0.2$. In Fig. \ref{fig_dumbell_results}(e), the changing direction by the motion by mean curvature is indicated by arrows.
Figure \ref{fig_star_results} shows the evolution of the star shape created by Eq. \eqref{2Dstarini}. The parameters are used as mentioned at the beginning of this section and $T=0.0325$. As shown in Fig. \ref{fig_star_results}, the tips of the star move inward and the gaps between the tips move outward. When it changes to the shape of a circle, the change of the radius can be predicted as shown in Fig. \ref{fig_circle_results}.
Figure \ref{fig_torus_results} shows the evolution of the torus shape of Eq. \eqref{2Dtorusini} with $T=0.0575$, $R_1=0.4$ and $R_2=0.3$. Because the inner circle has a larger curvature, it shrinks faster than the outer circle, and after the inner circle disappears, the change of radius over time can be measured as shown in Fig.\ref{fig_circle_results}.
Figure \ref{fig_complex_results} shows the evolution of a maze shape. We use $N_x=N_y=100$, $\epsilon_5$, and $T=0.04$. As shown in Fig. \ref{fig_complex_results}, we obtain the results of shrinking while maintaining its initial shape.
\clearpage
The last simulation in 2D is phase separation with a random initial condition Eq. \eqref{2Dtrnadini}. In Fig. \ref{fig_separation_results}, starting with random values with 0.1 amplitude, but over time, phase separation occurs with values $-1$ to $1$.

According to the results in Table \ref{tab:1}, the speed gap between CPU and GPU is significant. In 2D, that is up to 251.6 times the difference between the Python:CPU and the Pytorch:GPU codes. Also, Pytorch:GPU tensors make the model up to 4.73 times faster than Pytorch:CPU tensors in the same code. 
\begin{table}[!h]
 \begin{center}
   \caption{Runtime Result in 2D(sec). The values in parentheses describe how many times the difference is based on GPU:Pytorch time for each test.}\label{tab:1}
   \resizebox{\textwidth}{!}{
   \begin{tabular}{c c c c c c c}
   \hline
       &  \multicolumn{6}{c}{Initial value} \\
     \cline{2-7}
     $$           & $circle$ &  $dumbbell$ & $star$ & $torus$ & $maze$ & $separation$\\ 
     $Iterations$ & 12001    &  15001     & 13001  & 23001   & 4001   & 12001\\ 
     \hline
     CPU:Python  &  680.09(185.31)  & 1736.05(251.6)    & 729.74(128.47) & 1326.82(165.65) & 60.15(52.3) & 674.63(119.62)\\
     CPU:Pytorch &  14.73(4.01)   & 30.67(4.44)      & 17.60(3.1)  & 37.85(4.73)   & 2.98(2.59)  & 17.79(3.15) \\
     \hline 
     GPU:Pytorch &  3.67    & 6.90       & 5.68   & 8.01    & 1.15  & 5.64\\
     \hline
   \end{tabular}}
   \vspace{1ex}
 \end{center}
\end{table} 

\subsection{Simulations in 3 dimensional space} 
Three dimensional simulations can be considered as an extension of two dimensions tests. Unless otherwise stated, we use the following parameters: mesh size $N_x=N_y=N_z=100$, space step size $h=1/N_x$, time step size $\Delta t = 0.1h^2$, and computational domain $\Omega = (0,1)\times(0,1)\times(0,1)$. As in the 2D simulation, we find an appropriate $\epsilon_m$ by comparing the numerical solution with the exact solution (using in Eq. \eqref{exactradiusEQ}) for the radius of the sphere initial condition in Eq. \eqref{3Dsphereini}. 

\begin{figure}[htbp!]
\centering
\includegraphics[width=0.85\columnwidth]{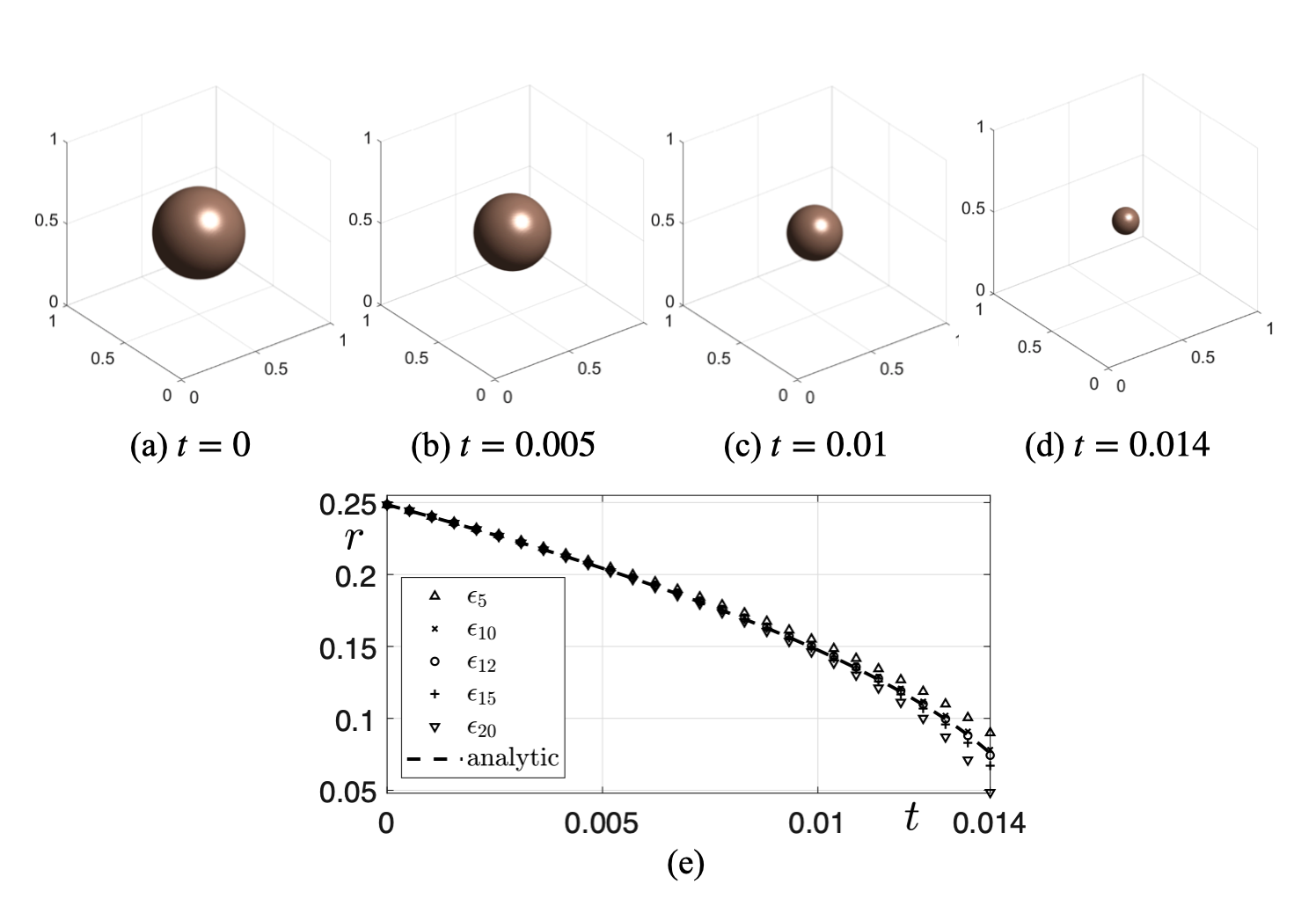}
\caption{(a)-(d) The changes in the circle over time are shown, and each time is described below in the figures. (e) Change of radius for various $\epsilon_m$.} \label{fig_sphere_results} 
\end{figure}

Figure \ref{fig_sphere_results}(a)-(d) show the time evolution results for $\epsilon_{12}$, and (e) presents the comparison of the numerical and exact solutions of radius for various epsilons. As shown in Fig. \ref{fig_sphere_results}, when $\epsilon_{12}$ is used, it is found that the exact solution and the numerical solution are most similar. Therefore, the other tests are performed using $\epsilon_{12}$.

\begin{figure}[htbp!]
\centering
\includegraphics[width=0.85\columnwidth]{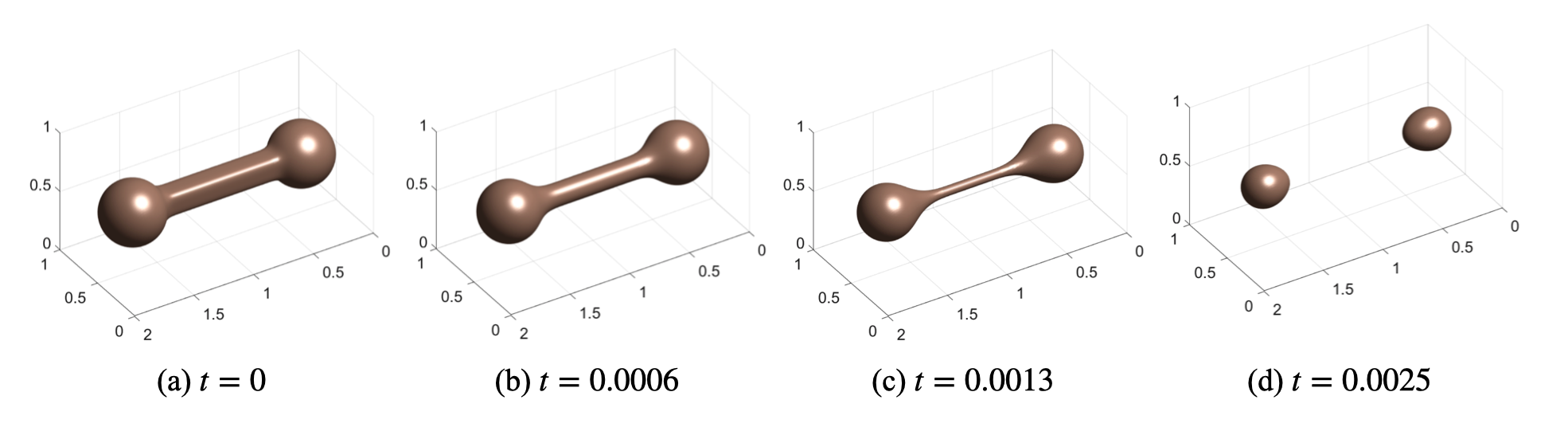}
\caption{(a)-(d) Time evolution of dumbbell shape. The handle shrinks quickly and breaks occurs, because the curvature of the handle is larger than that of both spheres} \label{fig_3Ddumbell_results} 
\end{figure}

\begin{figure}[htbp!]
\centering
\includegraphics[width=0.85\columnwidth]{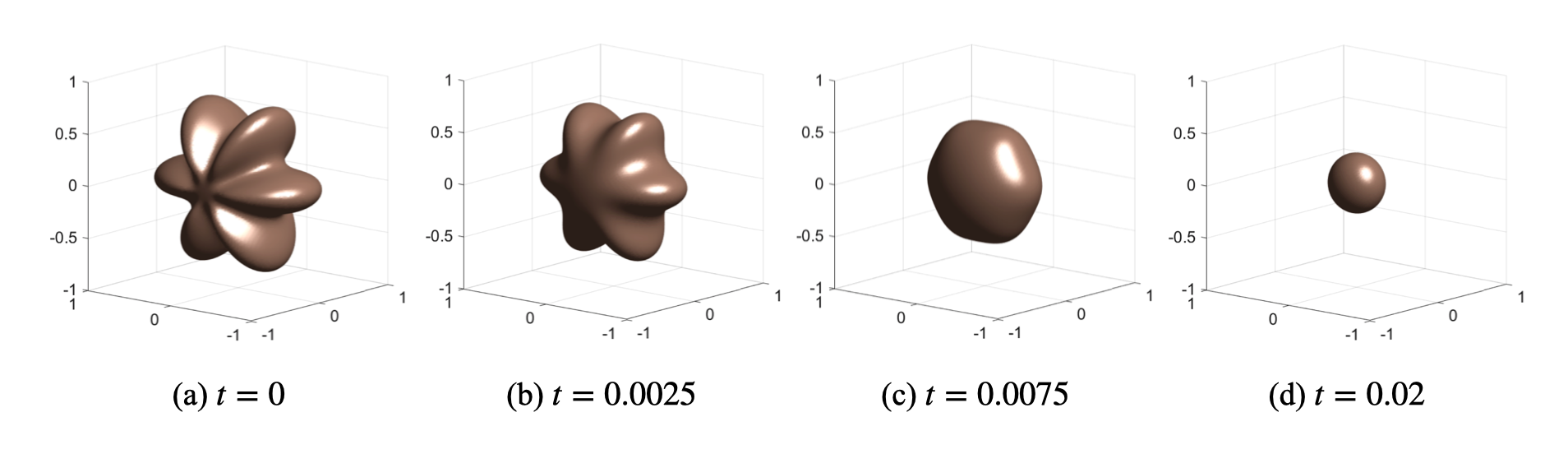}
\caption{Time evolution of a star shape. The shape changes with the mean curvature flow.} \label{fig_3Dstar_results} 
\end{figure}

\begin{figure}[htbp!]
\centering
\includegraphics[width=0.85\columnwidth]{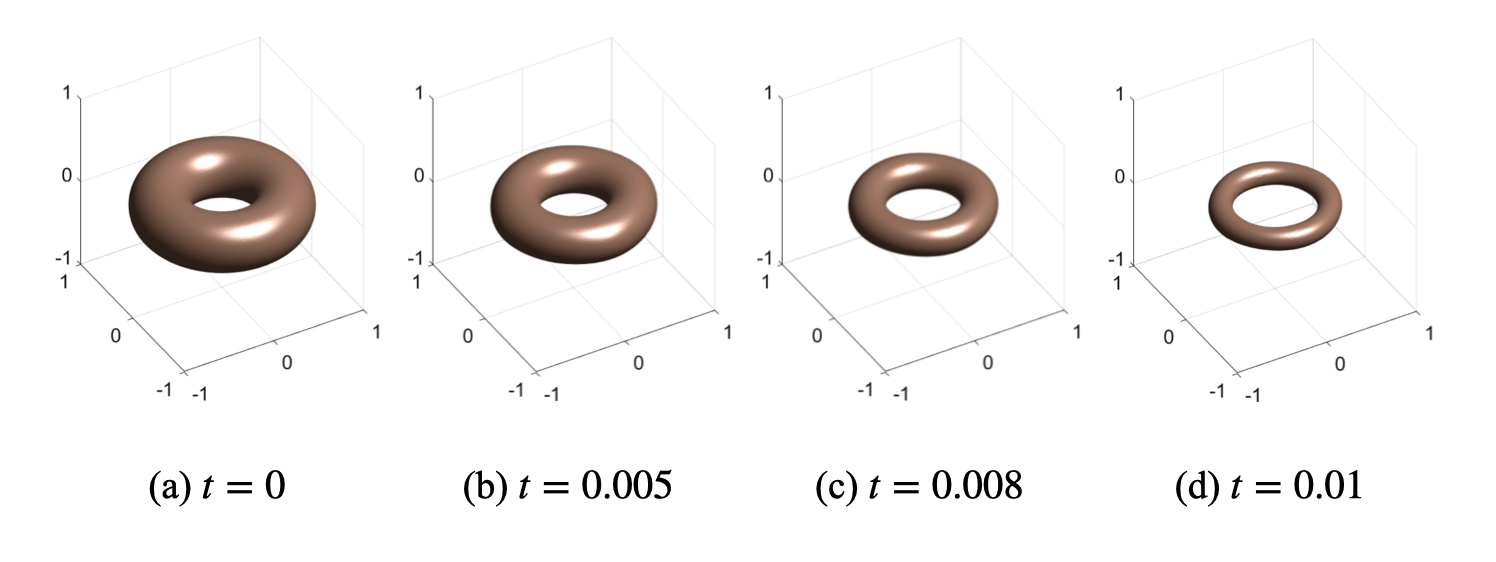}
\caption{Time evolution of a torus shape. Contrary to the 2D results, the inner circle increases due to the effect of the mean curvature.} \label{fig_3Dtorus_results} 
\end{figure}

\begin{figure}[htbp!]
\centering
\includegraphics[width=0.85\columnwidth]{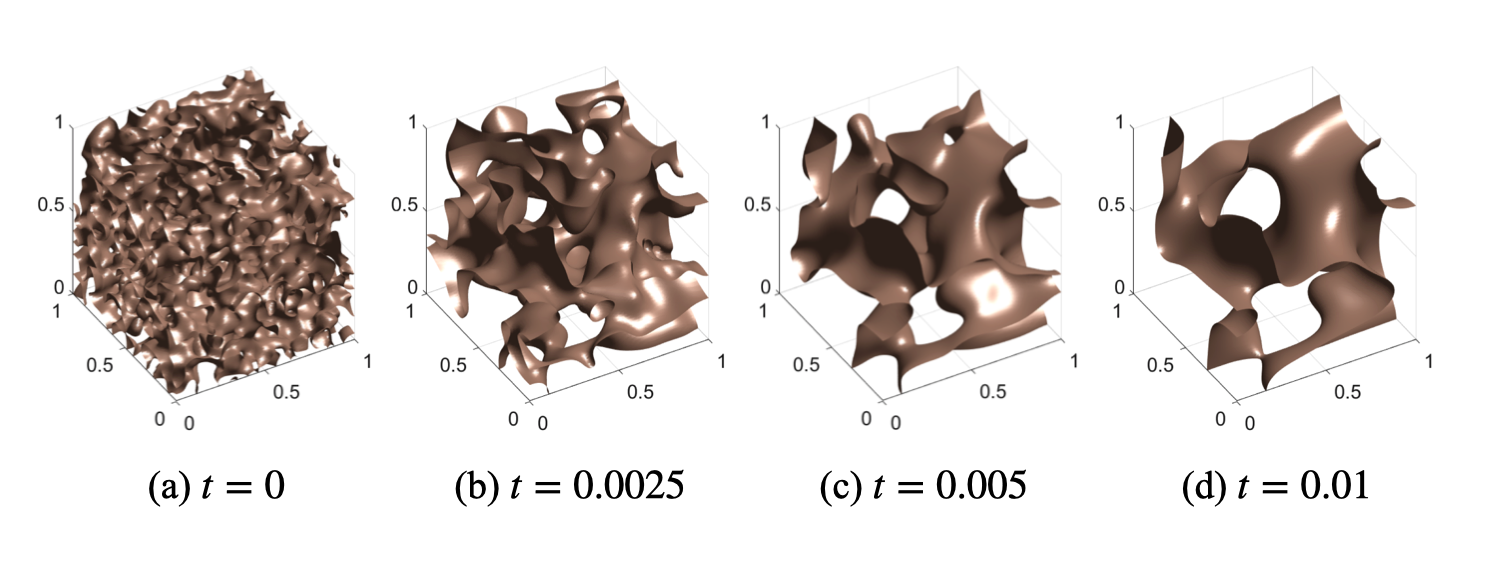}
\caption{Time evolution of phase separation with a random initial condition.} \label{fig_3Dseparation_results} 
\end{figure}

Figure \ref{fig_3Ddumbell_results} shows the temporal evolution of the dumbbell shape and the initial condition is shown in Eq. \eqref{3Ddumbellini}. For resolution, we use $N_x=200$ and $N_y=N_z=100$ on the computational domain $\Omega=(0,2)\times(0,1)\times(0,1)$. The final time $T$ of the simulation is $T=0.0025$ and $R_0=0.25$. The tendency of the 3D dumbbell motion is different from the result of 2D dumbbell. The reason is that in the initial shape in Fig. \ref{fig_3Ddumbell_results}(a), the radius of the handle is much smaller than the radius of the spheres at both ends. Therefore, the handle part shrinks faster than end of spheres and breaks as in Fig. \ref{fig_3Ddumbell_results}.
Figure \ref{fig_3Dstar_results} shows the evolution of the star shape shown in Eq. \eqref{3Dstarini} on the computational domain $
\Omega=(-1,1)\times(-1,1)\times(-1,1)$. The parameters are used as mentioned at the beginning of this section and $T=0.02$. As shown in Fig. \ref{fig_3Dstar_results}, similar to the 2D result, the tips of the star move inward and the gaps between the tips move outward.
Figure \ref{fig_3Dtorus_results} shows the evolution of the torus shape with initial condition in Eq. \eqref{3Dtorusini} using $R_1=0.3$ and $R_2=0.3$ and $T=0.01$ on the computational domain $
\Omega=(-1,1)\times(-1,1)\times(-1,1)$. Contrary to the 2D result, the inner circle becomes larger because the radius of the minor circle exists. Therefore, the radius of the major circle and the minor circle have different each curvature. Since the radius of the minor circle is smaller than major circle, the mean curvature drives the motion into the inside of the torus as shown in Fig. \ref{fig_3Dtorus_results}.
Figure \ref{fig_3Dcomplex_results} shows the evolution of a maze shape. The initial condition is described in Appendix code. In this simulation, we use $T=0.0175$ and computational domain $\Omega=(-1,1)\times(-1,1)\times(-1,1)$. As shown in Fig. \ref{fig_3Dcomplex_results}, we obtain the results of shrinking while maintaining its initial shape in Fig. \ref{fig_3Dcomplex_results}(b), then merges and shrinks in Fig. \ref{fig_3Dcomplex_results}(c)-(d). If we use different $\epsilon_m$, it can shrink while preserving its initial shape.
\begin{figure}[htbp!]
\centering
\includegraphics[width=0.85\columnwidth]{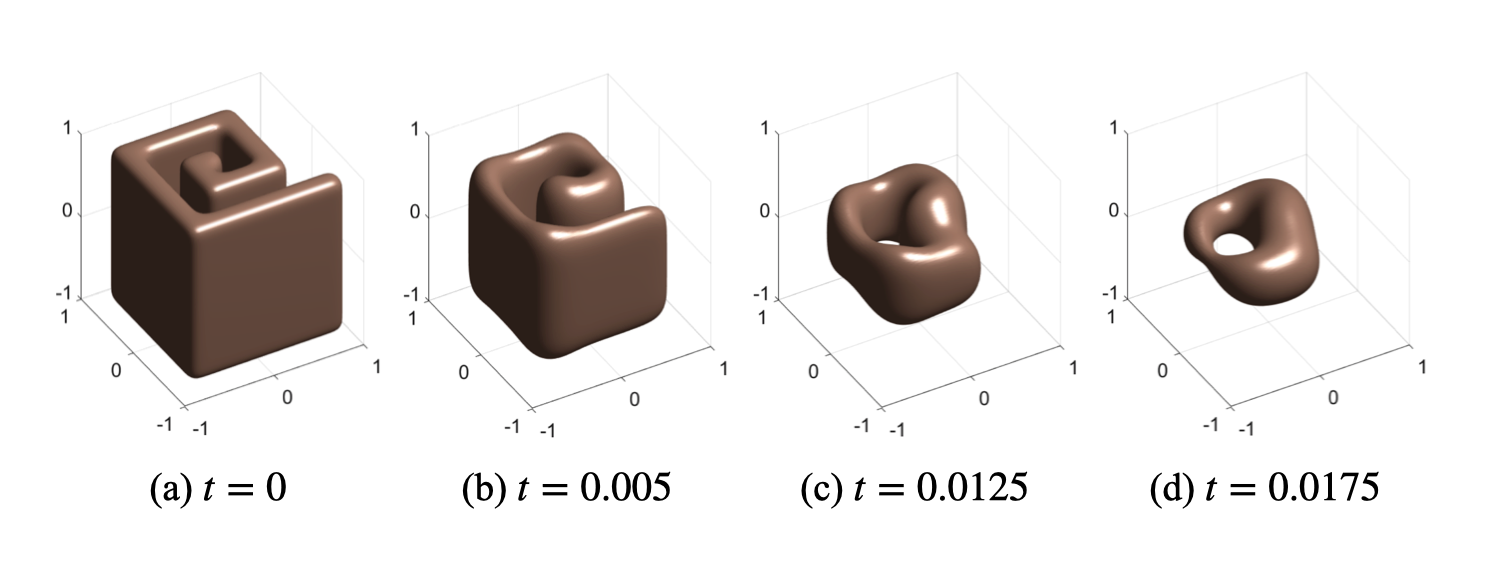}
\caption{Time evolution of a maze shape on $\Omega=(-1,1)\times(-1,1)\times(-1,1)$, and $T=0.0175$.} \label{fig_3Dcomplex_results} 
\end{figure}
The last simulation in 3D is phase separation with random initial condition in Eq. \eqref{3Drandini}. In Fig. \ref{fig_3Dseparation_results}, starting with random values with 0.1 amplitude, but over time, phase separation occurs with values $-1$ to $1$.

According to the results in Table \ref{tab:2}, the speed gap between CPU and GPU operations is significant. Of course, it will be faster by performing GPU calculations, but we proposed a structure using padding and convolution operation in performing GPU calculations on AC equations. And the results are demonstrated by verifying the results (Tab. \ref{tab:3}) with Python code which has the same mathematically meaning. In 3D, the GPU performance is much more overwhelming than the CPUs. For instance, the GPU code is 4766 times faster than the Python code in the torus problem. Also, GPU tensors make the model up to 76 times faster than CPU tensors in the same code.  The values in parentheses describe how many times the difference is based on GPU:Pytorch time for each test.

\begin{table}[htbp]
 \begin{center}
   \caption{Runtime Result in 3D(sec). The values in parentheses describe how many times the difference is based on GPU:Pytorch time for each test.}\label{tab:2}
   \resizebox{\textwidth}{!}{
   \begin{tabular}{c c c c c c c}
   \hline
       &  \multicolumn{6}{c}{Initial value} \\
     \cline{2-7}
     $ $          & $sphere$ &  $dumbbell$ & $star$ & $torus$  & $maze$ & $separation$\\ 
     $Iterations$ &  2001    &  2001      & 2001   & 1201     & 2401   & 2001\\ 
     \hline
     CPU:Python  &  4022.38(3944) & 8050.45(4087)    & 4014.83(4015) & 2478.22(4766) & 4844.24(3814) & 4163.17(4042)\\
     CPU:Pytorch &  65.95(65)   & 79.23(40)      & 65.22(65)   & 39.56(76)   & 78.72(62)   & 65.80(64) \\
     \hline 
     GPU:Pytorch &  1.02    & 1.97       & 1.00   & 0.52     & 1.27    & 1.03\\
     \hline
   \end{tabular}}
   \vspace{1ex}
 \end{center}
\end{table} 

\section{Conclusions} \label{con}

In this paper, we proposed a structure using padding and convolution operation for the GPU calculation of the Allen--Cahn equation. We increased the simulation speed and demonstrated the validity of the proposed structure by verifying the result with Python code that has the same mathematically meaning. We solved the Allen--Cahn equation by the explicit finite difference method and compared the runtime results between CPU and GPU algorithms. The errors of CPU:Python and GPU:Pytorch are less than $1.0e$-$6$ for the given initial conditions. Also, we showed that our GPU code is up to 251.6 and 4765.81 times faster than the CPU codes in 2D and 3D, respectively.  By showing these results, accuracy and efficiency have been demonstrated. Various numerical simulations were presented to confirm that the Allen--Cahn equation follows the motion by mean curvature and phase separation in 2D and 3D space. In this way, we can build a fast algorithm for any differential equations using the finite difference method by efficient programmatic code structures.


\end{document}